\documentclass[draft]{amsart}

\usepackage{amsthm}
\usepackage[leqno]{amsmath}
\usepackage{latexsym,amsfonts,amssymb}
\usepackage{hyperref}

\newcommand{\numberseries}{\mdseries}   

\newlength{\thmtopspace}                
\newlength{\thmbotspace}                
\newlength{\thmheadspace}               
\newlength{\thmindent}                  

\newtheoremstyle{bfupright head,slanted body}
                {\thmtopspace}{\thmbotspace}
                {\slshape}{\thmindent}{\bfseries}{.}{\thmheadspace}
                {{\numberseries \thmnumber{(#2) }}\thmnote{#3}}

\newtheoremstyle{bfupright head,upright body}
                {\thmtopspace}{\thmbotspace}
                {\upshape}{\thmindent}{\bfseries}{.}{\thmheadspace}
                {{\numberseries \thmnumber{(#2) }}\thmnote{#3}}

\newtheoremstyle{fixed bf head,upright body}
                {\thmtopspace}{\thmbotspace}{\upshape}
                {\thmindent}{\bfseries}{.}{\thmheadspace}
                {{\numberseries \thmnumber{(#2) }}\thmname{#1}\thmnote{ (#3)}}

\theoremstyle{bfupright head,slanted body}
\newtheorem{res}{}[section]             \newtheorem*{res*}{}

\theoremstyle{bfupright head,upright body}
\newtheorem{bfhpg}[res]{}               \newtheorem*{bfhpg*}{}

\theoremstyle{fixed bf head,upright body}
\newtheorem{rmk}[res]{Remark}           \newtheorem*{rmk*}{Remark}

\newcommand{\pgref}[1]{(\ref{#1})}
\newcommand{\secref}[2][Section~]{#1\ref{sec:#2}}
\renewcommand{\eqref}[1]{\pgref{eq:#1}}
\numberwithin{equation}{section}
 \newcommand{\e}{\varepsilon}
 \newcommand{\Koszul}[1][R]{\operatorname{K}^{#1}}
 \newcommand{\cls}[1]{\bar{#1}}
 \newcommand{\clC}[1]{\mathbf{C}(#1)}
 \newcommand{\clB}{\mathbf{B}}
 \newcommand{\clG}[1]{\mathbf{G}(#1)}

\def\urltilda{\kern -.15em\lower .7ex\hbox{\~{}}\kern .04em} 
\newcommand{\set}[2][\mspace{1mu}]{\{#1 #2 #1\}}
\newcommand{\kk}{\Bbbk}
\newcommand{\NN}{\mathbb{N}}
\newcommand{\qtext}[1]{\quad\text{#1}\quad}
\newcommand{\qqtext}[1]{\qquad\text{#1}\qquad}
\newcommand{\qand}{\qtext{and}}
\newcommand{\qqand}{\qqtext{and}}
\newcommand{\m}{\mathfrak{m}}
\newcommand{\is}{\cong}
\newcommand{\onto}{\twoheadrightarrow}
\newcommand{\lra}{\longrightarrow}
\newcommand{\pows}[2][k]{#1[\mspace{-2.3mu}[#2]\mspace{-2.3mu}]}
\newcommand{\Rhat}{\widehat{R}}
\newcommand{\mapdef}[4][\rightarrow]{\nobreak{#2\colon #3 #1 #4}}
\renewcommand{\H}[2][]{\operatorname{H}_{#1}(#2)}
\newcommand{\rnk}[2][k]{\operatorname{rank}_{#1}#2}
\newcommand{\tp}[3][R]{\nobreak{#2\otimes_{#1}#3}}
\newcommand{\Tor}[4][R]{\operatorname{Tor}^{#1}_{#2}(#3,#4)}

\begin{document}

\title[Local rings of embedding codepth $3$. Examples]{Local rings of
  embedding codepth $\mathbf{3}$. Examples}

\author[L.\,W. Christensen]{Lars Winther Christensen}

\address{Texas Tech University, Lubbock, TX 79409, U.S.A.}

\email{lars.w.christensen@ttu.edu}

\urladdr{http://www.math.ttu.edu/\urltilda lchriste}

\author[O. Veliche]{Oana Veliche}

\address{Northeastern University, Boston, MA~02115, U.S.A.}

\email{o.veliche@neu.edu}

\urladdr{http://www.math.northeastern.edu/\urltilda veliche}

\thanks{This research was partly supported by NSA grant H98230-11-0214
  (L.W.C.).}

\date{14 November 2012}





\begin{abstract}
  A complete local ring of embedding codepth $3$ has a minimal free
  resolution of length $3$ over a regular local ring. Such resolutions
  carry a differential graded algebra structure, based on which one
  can classify local rings of embedding codepth $3$. We give examples
  of algebra structures that have been conjectured not to occur.
\end{abstract}

\maketitle

\thispagestyle{empty}

\section{Introduction} 

\noindent A classification of commutative noetherian local rings of
embedding codepth $c \le 3$ took off more than a quarter of a century
ago. Up to completion, a local ring of embedding codepth $c$ is a
quotient of regular local ring $Q$ by an ideal $I$ of grade $c$, and
the classification is based on an algebra structure on
$\Tor[Q]{*}{Q/I}{k}$, where $k$ is the residue field of $Q$. The
possible isomorphism classes of these algebras were identified by
Weyman~\cite{JWm89} and by Avramov, Kustin, and Miller
\cite{AKM-88}. Significant restrictions on the invariants that
describe these isomorphism classes were recently worked out by
Avramov~\cite{LLA}. Here is a pr\'ecis that will suffice for our
purposes.

Let $R$ be a commutative noetherian local ring with maximal ideal $\m$
and residue field $k=R/\m$. Denote by $e$ the minimal number of
generators of $\m$ and by $d$ the depth of $R$. The number $e$ is
called the \emph{embedding dimension} of $R$, and $c=e-d$ is the
\emph{embedding codepth}. By Cohen's Structure Theorem the $\m$-adic
completion $\Rhat$ of $R$ has the form $\Rhat = Q/I$, where $Q$ is a
complete regular local ring with the same embedding dimension and
residue field as $R$; we refer to $I$ as the \emph{Cohen ideal} of
$R$.

The projective dimension of $\Rhat$ over $Q$ is $c$, by the
Auslander--Buchsbaum Formula. From now on let $c\le 3$; the minimal
free resolution $F$ of $\Rhat$ over $Q$ then carries a differential
graded algebra structure. It induces a graded algebra structure on
$\tp[Q]{F}{k} = \H{\tp[Q]{F}{k}} = \Tor[Q]{*}{\Rhat}{k}$, which
identifies $R$ as belonging to one of six (parametrized) classes,
three of which are called $\clB$, $\clC{c}$, and $\clG{r}$ for $r\ge
2$.  The ring $R$ is in $\clC{c}$ if and only if it is an embedding
codepth $c$ complete intersection. If $R$ is Gorenstein but not a
complete intersection, then it belongs to the class $\clG{r}$ with $r
=\mu(I)$, the minimal number of generators of the Cohen ideal. Work of
J.~Watanabe~\cite{JWt73} shows that $\mu(I)$ is odd and at least
$5$. Brown~\cite{AEB87} identified rings in $\clB$ of type $2$, and
thus far no other examples of $\clB$ rings have been known. Rings in
$\clG{r}$ that are not Gorenstein---rings in $\clG{3}$ and $\clG{2n}$
for $n\in\NN$ in particular---have also been elusive; in fact, it has
been conjectured \cite[3.10]{LLA} that every ring in $\clG{r}$ would
be Gorenstein and, by implication, that the classes $\clG{3}$ and
$\clG{2n}$ would be empty.

In this paper we provide examples of some of the sorts of rings that
have hitherto dodged detection; the precise statements follow in
Theorems I and II below.

\begin{res*}[Theorem I]
  Let $\kk$ be a field, set $Q=\pows[\kk]{x,y,z}$, and consider these
  ideals in~$Q$:
  \begin{align*}
    \mathfrak{g}_1 &= (xy^2,xyz,yz^2,x^4-y^3z, xz^3-y^4)\\
    \mathfrak{g}_2 &= \mathfrak{g}_1 + (x^3y-z^4)\\
    \mathfrak{g}_3 &= \mathfrak{g}_2 + (x^2z^2)\\
    \mathfrak{g}_4 &= \mathfrak{g}_3 + (x^3z)\;.
  \end{align*}
  Each algebra $Q/\mathfrak{g}_n$ has embedding codepth 3 and type
  $2$, and $Q/\mathfrak{g}_n$ is in $\clG{n+1}$.
\end{res*}

The theorem provides counterexamples to the conjecture mentioned
above: The classes $\clG{2}$, $\clG{3}$, and $\clG{4}$ are not empty,
and rings in $\clG{5}$ need not be Gorenstein.  The second theorem
provides examples of rings in $\clB$ of type different from~$2$.

\begin{res*}[Theorem II]
  Let $\kk$ be a field, set $Q=\pows[\kk]{x,y,z}$, and consider these
  ideals in~$Q$:
  \begin{align*}
    \mathfrak{b}_1 &= (x^3,x^2y,yz^2,z^3)\\
    \mathfrak{b}_2 &= \mathfrak{b}_1 + (xyz)\\
    \mathfrak{b}_3 &= \mathfrak{b}_2 + (xy^2-y^3)\\
    \mathfrak{b}_4 &= \mathfrak{b}_3 + (y^2z)\;.
  \end{align*}
  Each algebra $Q/\mathfrak{b}_n$ has embedding codepth $3$ and
  belongs to $\clB$. The algebras $Q/\mathfrak{b}_1$ and
  $Q/\mathfrak{b}_2$ have type $1$ while $Q/\mathfrak{b}_3$ and
  $Q/\mathfrak{b}_4$ have type $3$.
\end{res*}

The algebras $Q/\mathfrak{b}_1$ and $Q/\mathfrak{b}_2$ have embedding
codimension $2$, which is the largest possible value for a
non-Gorenstein ring of embedding codepth $3$ and type $1$. The
algebras $Q/\mathfrak{b}_3$ and $Q/\mathfrak{b}_4$ are artinian; that
is, they have embedding codimension~$3$. An artinian local ring of
codepth $3$ and type $2$ belongs to $\clB$ only if the minimal
number of generators of its Cohen ideal is odd and at least $5$; see
\cite[3.4]{LLA}. In that light it appears noteworthy that
$\mathfrak{b}_3$ and $\mathfrak{b}_4$ are minimally generated by $6$
and $7$ elements, respectively.
\begin{equation*}
  \ast \ \ast \ \ast
\end{equation*}
In preparation for the proofs, we recall a few definitions and
facts. Let $Q$ be a regular local ring with residue field $k$ and let
$I$ be an ideal in $Q$ of grade $3$. The quotient ring $R=Q/I$ has
codepth $3$ and its minimal free resolution over $Q$ has the form
\begin{equation*}
  F \; = \; 0 \lra Q^n \lra Q^{m} \lra Q^{l+1} \lra Q \lra 0\;,
\end{equation*}
where $n$ is the type of $R$ and one has $m=n+l$ and $l+1=\mu(I)$. It
has a structure of a graded-commutative differential graded algebra;
this was proved by Buchsbaum and Eisenbud \cite[1.3]{DABDEs77}. While
this structure is not unique, the induced graded-commutative algebra
structure on $A=\H{\tp[Q]{F}{k}}$ is unique up to isomorphism. Given
bases
\begin{equation}
  \label{eq:bases}
  \begin{aligned}
    \mathbf{e}_1,\ldots,\mathbf{e}_{l+1} &\quad\text{for $A_1$},\\
    \mathbf{f}_1,\ldots,\mathbf{f}_m &\quad\text{for $A_2$},\;\text{and}\\
    \mathbf{g}_1,\ldots,\mathbf{g}_n &\quad\text{for $A_3$}
  \end{aligned}
\end{equation}
graded-commutativity yields
\begin{equation}
  \label{eq:grcom}
  \begin{alignedat}{2}
    -\mathbf{e}_i\mathbf{e}_j &= \mathbf{e}_j\mathbf{e}_i \text{ and } \mathbf{e}_i^2=0& \quad
    &\text{for all } 1\le i,j\le l+1\;;\\
    \mathbf{e}_i\mathbf{f}_j &= \mathbf{f}_j\mathbf{e}_i && \text{for
      all } \ 1\le i\le l+1 \ \text{ and all } \ 1\le j\le m\;.
  \end{alignedat}
\end{equation}

We recall from \cite[2.1]{AKM-88} the definitions of the classes
$\clG{r}$ with $r\ge 2$\footnote{ One does not define $\clG{1}$
  because it would overlap with another class called
  $\mathbf{H}(0,1)$.} and $\clB$.

The ring $R$ belongs to $\clG{r}$ if there is a basis \eqref{bases}
for $A_{\ge 1}$ such that one has
\begin{equation}
  \label{eq:G}
  \mathbf{e}_i\mathbf{f}_i = \mathbf{g}_1 
  \qqtext{ for all $1 \le i \le r$}
\end{equation}
and all other products of basis elements not fixed by \eqref{G} via
\eqref{grcom} are zero.

The ring $R$ belongs to $\clB$ if there is a basis \eqref{bases} for
$A_{\ge 1}$ such that one has
\begin{equation}
  \label{eq:B}
  \mathbf{e}_1\mathbf{e}_2 = \mathbf{f}_3 \qqand 
  \begin{aligned}
    \mathbf{e}_1\mathbf{f}_1 &= \mathbf{g}_1\\
    \mathbf{e}_2\mathbf{f}_2 &= \mathbf{g}_1
  \end{aligned}
\end{equation}
and all other products of basis elements not fixed by \eqref{B} via
\eqref{grcom} are zero.

The proofs of both theorems use the fact that the graded algebra $A$ is
isomorphic to the Koszul homology algebra over $R$. We fix notation
for the Koszul complex. Let $R$ be of embedding dimension $3$ and let
$\m = (x,y,z)$ be its maximal ideal. We denote by $\Koszul$ the Koszul
complex over the canonical homomorphism
$\mapdef[\onto]{\pi}{R^3}{\m}$. It is the exterior algebra of the rank
$3$ free $R$-module with basis $\e_x$, $\e_y$, $\e_z$, endowed with
the differential induced by $\pi$. For brevity we set
\begin{equation*}
  \e_{xy} = \e_x\wedge\e_y\,,\quad\e_{xz} = \e_x\wedge\e_z\,,\quad 
  \e_{yz} = \e_y\wedge\e_z 
  \qand  \, \e_{xyz} = \e_x\wedge\e_y\wedge\e_z\;.
\end{equation*}
The differential is then given by,
\begin{alignat*}{3}
  \partial(\e_{x}) &= x \qquad &\partial(\e_{xy}) &= x\e_{y} - y\e_{x}\\
  \partial(\e_{y}) &= y & \partial(\e_{xz}) &= x\e_{z} -
  z\e_{x}&\qquad
  \partial(\e_{xyz}) &= x\e_{yz} - y\e_{xz} + z\e_{xy}\;,\\
  \partial(\e_{z}) &= z &\partial(\e_{yz}) &= y\e_{z} - z\e_{y}
\end{alignat*}
and it makes $\Koszul$ into a graded-commutative differential graded
algebra. The induced algebra structure on $\H{\Koszul}$ is
graded-commutative, and there is an isomorphism of graded algebras $A
\is \H{\Koszul}$; see \cite[(1.7.1)]{LLA}.

Note that the homology module $\H[3]{\Koszul}$ is isomorphic, as a
$k$-vector space, to the socle of $R$, that is, the ideal $(0:\m)$. To
be precise, if $s_1,\ldots,s_n$ is a basis for the socle of $R$, then
the homology classes of the cycles $s_1\e_{xyz},\ldots,s_n\e_{xyz}$ in
$\Koszul_3$ form a basis for $\H[3]{\Koszul}$. From the isomorphism $A
\is \H{\Koszul}$ one gets, in particular, $\sum_{i=0}^3(-1)^i
\rnk[\kk]{\H[i]{\Koszul}} = 0$, as the ranks of the Koszul homology
modules equal the ranks of the free modules in $F$. To sum up one has,
\begin{equation}
  \label{eq:rank}
  \begin{aligned}
    \rnk[\kk]{\H[1]{\Koszul}} &= l+1 = \mu(I)\\
    \rnk[\kk]{\H[2]{\Koszul}} &= l + \rnk[\kk]{\H[3]{\Koszul}}\\
    \rnk[\kk]{\H[3]{\Koszul}} &= \rnk[\kk]{(0:\m)}\;.
  \end{aligned}
\end{equation}

Theorem I is proved in Sections \secref[]{I2}--\secref[]{I45} and
Theorem II in Sections \secref[]{J1}--\secref[]{J34}. For each
quotient algebra $Q/\mathfrak{g}_n$ and $Q/\mathfrak{b}_n$ we shall
verify that the Koszul homology algebra has the desired multiplicative
structure as described in \eqref{G} and \eqref{B}, and we shall
determine the type of the quotient algebra. As described above, the
latter means determining the socle rank, as each of these algebras has
depth $0$.

\section{ Proof that $Q/\mathfrak{g}_1$ is a type $2$ algebra in
  $\clG{2}$}
\label{sec:I2}

\noindent
The ideal $\mathfrak{g}_1=(xy^2,xyz,yz^2,x^4-y^3z, xz^3-y^4)$ in
$Q=\pows[\kk]{x,y,z}$ is generated by homogeneous elements, so
$R=Q/\mathfrak{g}_1$ is a graded $\kk$-algebra.  For $n\ge 0$ denote
by $R_n$ the subspace of $R$ of homogeneous polynomials of degree
$n$. It is simple to verify that the elements in the second column
below form bases for the subspaces $R_n$; for convenience, the third
column lists the relations among non-zero monomials.
\begin{equation}
  \label{eq:I2}
  \text{\begin{tabular}{lll}
      $R_0$ & $1$\\
      $R_1$ & $x,\, y,\, z$\\
      $R_2$ & $x^2,\, xy,\, xz,\, y^2,\, yz,\, z^2$\\
      $R_3$ & $x^3,\, x^2y,\, x^2z,\, xz^2,\, y^3,\, y^2z,\, z^3\quad$\\
      $R_4$ & $x^4,\, x^3y,\, x^3z,\, x^2z^2,\, xz^3,\, z^4$& $y^4=xz^3,\, y^3z=x^4$ \\
      $R_5$ & $x^4y,\, x^3z^2,\, z^5$& $xz^4=y^4z=x^4y$\\
      $R_{n\ge 6}$ & $z^n$\\
    \end{tabular}}
\end{equation}
Set $A=\H{\Koszul}$; we shall verify that $A$ has the multiplicative
structure described in \eqref{G} with $r=2$, and that $R$ has socle
rank $2$.

\begin{bfhpg*}[A basis for $A_3$]
  From \eqref{I2} it is straightforward to verify that the socle of
  $R$ is generated by $x^4y$ and $x^3z^2$, so it has rank 2 and the
  homology classes of the cycles
  \begin{equation}
    \label{eq:I2g}
    g_1 = x^4y\e_{xyz} \qqand g_2 = x^3z^2\e_{xyz}
  \end{equation}
  form a basis for $A_3$. As there are no non-zero boundaries in
  $\Koszul_3$, the homology classes $\cls{g}_1$ and $\cls{g}_2$
  contain only $g_1$ and $g_2$, respectively, and the bar merely
  signals that we consider the cycles as elements in $A$ rather than
  $\Koszul$.
\end{bfhpg*}

As the ideal $\mathfrak{g}_1$ is minimally generated by $5$ elements,
one has $\rnk[\kk]{A_1}=5$, and hence $\rnk[\kk]{A_2}=6$ by
\eqref{rank}.

\begin{bfhpg*}[A basis for $A_2$]
  It is elementary to verify that the next elements in $\Koszul_2$ are
  cycles.
  \begin{equation}
    \label{eq:I2f}
    \begin{aligned}
      f_1 &= -yz\e_{xy} + y^2\e_{xz}\\
      f_2 &= yz\e_{xz}\\
      f_3 &= x^3z\e_{xy}\\
      f_4 &= xz^3\e_{xy}\\
      f_5 &= (x^3y-z^4)\e_{xy}\\
      f_6 &= x^3z^2\e_{xz}
    \end{aligned}
  \end{equation}
  To see that their homology classes form a basis for $A_2$ it
  suffices, as $A_2$ is a $\kk$-vector space of rank $6$, to verify
  that they are linearly independent modulo boundaries. A boundary in
  $\Koszul_2$ has the form
  \begin{equation}
    \label{eq:B2}
    \partial(a\e_{xyz}) = ax\e_{yz} - ay\e_{xz} + az\e_{xy}\;,
  \end{equation}
  for some $a$ in $R$.  As the differential is graded, one needs to
  verify that $f_1$ and $f_2$ are linearly independent modulo
  boundaries, that $f_3$, $f_4$, and $f_5$ are linearly independent
  modulo boundaries, and that $f_6$ is not a boundary.

  If $u$ and $v$ are elements in $\kk$ such that $uf_1+vf_2 =
  -uyz\e_{xy} + (uy^2+vyz)\e_{xz}$ is a boundary, then it has the form
  \eqref{B2} for some $a$ in $R_1$. In particular, one has $ax=0$, and
  that forces $a=0$; see \eqref{I2}. With this one has $uyz = 0$ and
  $uy^2+vyz=0$, whence $u=0=v$. Thus, $f_1$ and $f_2$ are linearly
  independent modulo boundaries.

  If $u$, $v$, and $w$ in $\kk$ are such that $uf_3+vf_4+wf_5 = (wx^3y
  + ux^3z + vxz^3 -wz^4)\e_{xy}$ is a boundary, then it has the form
  \eqref{B2} for some element
  \begin{equation*}
    a=a_1x^3 +a_2x^2y+a_3x^2z +a_4xz^2+ a_5y^3+a_6y^2z+a_7z^3
  \end{equation*}
  in $R_3$. From $az=wx^3y + ux^3z + vxz^3 - wz^4$ one gets $w=0$,
  $u=a_1$, and $v=a_4$; see~\eqref{I2}. From $ax=0$ one gets
  $a_1=0=a_4$, that is, $u=0=v$.

  Finally, $f_6$ is not a boundary as no element $a$ in $R_4$
  satisfies $-ay=x^3z^2$.
\end{bfhpg*}

\begin{bfhpg*}[A basis for $A_1$]
  The following elements are cycles in $\Koszul_1$.
  \begin{equation}
    \label{eq:I2e}
    \begin{aligned}
      e_1 &= x^3\e_x - y^2z\e_y\\
      e_2 &= z^3\e_x - y^3\e_y\\
      e_3 &= yz\e_x\\
      e_4 &= z^2\e_y\\
      e_5 &= y^2\e_x
    \end{aligned}
  \end{equation}
  The vector space $A_1$ has rank $5$, so as above the task is to show
  that $e_1,\ldots,e_5$ are linearly independent modulo boundaries. A
  boundary in $\Koszul_1$ has the form
  \begin{equation}
    \label{eq:B1}
    \partial(a\e_{xy} + b\e_{xz} + c\e_{yz}) = -(ay+bz)\e_{x} 
    + (ax-cz)\e_{y} + (bx+cy)\e_{z}\;,
  \end{equation}
  for $a$, $b$, and $c$ in $R$. As above, the fact that the
  differential is graded allows us to treat cycles with coefficients
  in $R_2$ and $R_3$ independently.

  If $u$ and $v$ are elements in $\kk$ such that $ue_1 + ve_2 =
  (ux^3+vz^3)\e_x - (uy^2z+vy^3)\e_y$ is a boundary in $\Koszul_1$,
  then it has the form \eqref{B1} for elements
  \begin{align*}
    a &= a_1x^2 + a_2xy + a_3xz + a_4y^2 + a_5yz + a_6z^2\;,\\
    b &= b_1x^2 + b_2xy + b_3xz + b_4y^2 + b_5yz + b_6z^2\;,\qand\\
    c &= c_1x^2 + c_2xy + c_3xz + c_4y^2 + c_5yz + c_6z^2
  \end{align*}
  in $R_2$. The equality $cz-ax=uy^2z + vy^3$ yields $v=0$ and
  $u=c_4$, and from $bx+cy=0$ one gets $c_4=0$.

  If $u$, $v$, and $w$ are elements in $\kk$ such that $ue_3 + ve_4 +
  we_5 = (uyz+wy^2)\e_x + vz^2\e_y$ is a boundary, then it has the
  form \eqref{B1} for $a=a_1x+a_2y+a_3z$, $b=b_1x+b_2y+b_3z$, and
  $c=c_1x+c_2y+c_3z$ in $R_1$. This yields equations:
  \begin{align*}
    -a_1xy - a_2y^2 - (a_3+b_2)yz - b_1xz - b_3z^2 &= uyz+wy^2\;,\\
    a_1x^2 + a_2xy + (a_3-c_1)xz - c_2yz - c_3z^2 &= vz^2\;,\qand\\
    b_1x^2 + (b_2+c_1)xy + b_3xz + c_2y^2 + c_3yz &= 0\;.
  \end{align*}
  From the last equation one gets, in particular, $c_3=0$ and
  $b_2+c_1=0$. The second one now yields $v=0$, $a_2=0$, and
  $a_3=c_1$. With these equalities, the first equation yields $w=0$
  and $u=0$.
\end{bfhpg*}

\begin{bfhpg*}[The product $A_1 \cdot A_1$]
  To determine the multiplication table $A_1\times A_1$ it is by
  \eqref{grcom} sufficient to compute the products $e_ie_j$ for $1\le
  i<j\le5$.  The product $e_3e_5$ is zero by graded-commutativity of
  $\Koszul$, and the following products are zero because the
  coefficients vanish in $R$; cf.~\eqref{I2}.
  \begin{alignat*}{2}
    e_1e_2 &= (-x^3y^3+y^2z^4)\e_{xy}&\qquad
    e_3e_4 &= yz^3\e_{xy}\\
    e_1e_3 &= y^3z^2\e_{xy}&
    e_4e_5 &= -y^2z^2\e_{xy}\\
    e_2e_5 &= y^5\e_{xy}
  \end{alignat*}
  Finally, the computations
  \begin{alignat*}{2}
    e_1e_4 &= x^3z^2\e_{xy} = \partial(x^3z\e_{xyz})\;,\\
    e_1e_5 &= y^4z\e_{xy}= xz^4\e_{xy} = \partial(xz^3\e_{xyz})\;,\\
    e_2e_3 &= y^4z\e_{xy} = xz^4\e_{xy} = \partial(xz^3\e_{xyz})\;,\qand\\
    e_2e_4 &= z^5\e_{xy} = \partial((z^4-x^3y)\e_{xyz})
  \end{alignat*}
  show that also the products $\cls{e}_1\cls{e}_4$,
  $\cls{e}_1\cls{e}_5$, $\cls{e}_2\cls{e}_3$, and $\cls{e}_2\cls{e}_4$
  in homology are zero. Thus, one has $A_1\cdot A_1=0$.
\end{bfhpg*}

\begin{bfhpg*}[The product $A_1 \cdot A_2$]
  Among the products $e_if_j$ for $1\le i\le 5$ and $1\le j \le 6$
  several are zero by graded-commutativity of $\Koszul$:
  \begin{gather*}
    e_3f_j \ \text{ and } \ e_5f_j \ \text{ for } 1\le j\le 6\\
    e_1f_3\,,\; e_1f_4\,,\; e_1f_5\,,\; e_2f_3\,,\; e_2f_4\,,\;
    e_2f_5\,,\; e_4f_3\,,\; e_4f_4\,, \text{ and }\, e_4f_5\;.
  \end{gather*}
  The following products are zero because the coefficients vanish in
  $R$.
  \begin{alignat*}{2}
    e_1f_2 &= y^3z^2\e_{xyz}&\qquad
    e_4f_1 &= -y^2z^2\e_{xyz}\\
    e_1f_6 &= x^3y^2z^3\e_{xyz}&
    e_4f_2 &= -yz^3\e_{xyz}\\
    e_2f_1 &= y^5\e_{xyz}&
    e_4f_6 &= -x^3z^4\e_{xyz}\\
    e_2f_6 &= x^3y^3z^2\e_{xyz}
  \end{alignat*}
  This leaves two products to compute, namely $e_1f_1 = y^4z\e_{xyz} =
  e_2f_2$.
\end{bfhpg*}

The computations above show that in terms of the $\kk$-basis
$\cls{e}_1,\ldots,\cls{e}_5$, $\cls{f}_1,\ldots,\cls{f}_6$,
$\cls{g}_1,\cls{g}_2$ for $A_{\ge 1}$ the non-zero products of basis
vectors are
\begin{equation}
  \label{eq:I2G}
  \cls{e}_1\cls{f}_1 = \cls{e}_2\cls{f}_2 = \cls{g}_1\;,
\end{equation}
whence $R$ belongs to $\clG{2}$.

\section{Proof that $Q/\mathfrak{g}_2$ is a type $2$ algebra in
  $\clG{3}$}
\label{sec:I3}

\noindent Set $R=Q/\mathfrak{g}_2$; as one has
$\mathfrak{g}_2=\mathfrak{g}_1+(x^3y-z^4)$ it follows from \eqref{I2}
that the elements listed below form bases for the subspaces $R_n$. As
in \eqref{I2} the third column records the relations among non-zero
monomials.
\begin{equation}
  \label{eq:I3}
  \text{\begin{tabular}{lll}
      $R_0$ & $1$\\
      $R_1$ & $x,\, y,\, z$\\
      $R_2$ & $x^2,\, xy,\, xz,\, y^2,\, yz,\, z^2$\\
      $R_3$ & $x^3,\, x^2y,\, x^2z,\, xz^2,\, y^3,\, y^2z,\, z^3\quad$\\
      $R_4$ & $x^4,\, x^3y,\, x^3z,\, x^2z^2,\, xz^3$ 
      & $y^4=xz^3,\,y^3z=x^4,\,z^4=x^3y$ \\
      $R_5$ & $x^4y,\,x^3z^2$ & $xz^4=y^4z=x^4y$\\
    \end{tabular}}
\end{equation}
Set $A=\H{\Koszul}$; we shall verify that $A$ has the multiplicative
structure described in \eqref{G} with $r=3$, and that $R$ has socle
rank $2$.

The next remark will also be used in later sections; loosely speaking,
it allows us to recycle the computations from \secref{I2} in the
analysis of $A$.

\stepcounter{res}
\begin{rmk}
  \label{notation}
  Let $\mathfrak{a} \subseteq \mathfrak{b}$ be ideals in $Q$. The
  canonical epimorphism $Q/\mathfrak{a} \onto Q/\mathfrak{b}$ yields a
  morphism of complexes $\Koszul[Q/\mathfrak{a}] \to
  \Koszul[Q/\mathfrak{b}]$. It maps cycles to cycles and boundaries to
  boundaries. To be explicit, let $E =
  \set{1,\e_{x},\e_{y},\e_{z},\e_{xy},\e_{xz},\e_{yz},\e_{xyz}}$ be
  the standard basis for either Koszul complex; if $\sum_{\e\in E}
  (q_\e +\mathfrak{a}) \e$ is a cycle (boundary) in
  $\Koszul[Q/\mathfrak{a}]$, then $\sum_{\e\in E} (q_\e +
  \mathfrak{b}) \e$ is a cycle (boundary) in
  $\Koszul[Q/\mathfrak{b}]$. By habitual abuse of notation, we write
  $x$, $y$, and $z$ for the cosets of $x$, $y$, and $z$ in any
  quotient algebra of $Q$, and as such we make no notational
  distinction between an element in $\Koszul[Q/\mathfrak{a}]$ and its
  image in $\Koszul[Q/\mathfrak{b}]$.
\end{rmk}

\begin{bfhpg*}[A basis for $A_{\ge 1}$]
  From \eqref{I3} it is straightforward to verify that the socle of
  $R$ is $R_5$, so it has rank $2$, and the homology classes of the
  cycles $g_1$ and $g_2$ from \eqref{I2g} form a basis for $A_3$.  The
  ideal $\mathfrak{g}_2$ is minimally generated by $6$ elements, so
  one has $\rnk[\kk]{A_1}=6$ and $\rnk[\kk]{A_2}=7$; see
  \eqref{rank}. Proceeding as in \secref{I2} it is straightforward to
  verify that the homology classes $e_1,\ldots,e_5$ from \eqref{I2e}
  together with the homology class of the cycle
  \begin{equation*}
    e_6 = x^2y\e_{x}-z^3\e_z
  \end{equation*}
  form a basis for $A_1$. Similarly, one verifies that the homology
  classes of $f_1$, $f_2$, $f_4$, and $f_6$ from \eqref{I2f} together
  with those of the cycles
  \begin{equation*}
    \begin{aligned}
      f_3 &= x^3\e_{xy} - z^3\e_{xz} + y^3\e_{yz}\;,\\
      f_5 &= z^4\e_{xz} - y^3z\e_{yz}\;, \text{ and}\\
      f_7 &= x^2y\e_{xy}+z^3\e_{yz}
    \end{aligned}
  \end{equation*}
  make up a basis for $A_2$.
\end{bfhpg*}

\begin{bfhpg*}[The product $A_1 \cdot A_1$]
  It follows from the computations in the previous section that one
  has $\cls{e}_i\cls{e}_j=0$ for $1\le i,j \le 5$.  To complete the
  multiplication table $A_1\times A_1$ it is by \eqref{grcom}
  sufficient to compute the products $e_ie_6$ for $1\le i \le 5$.  The
  products $e_1e_6$ and $e_2e_6$ are zero as one has $R_3\cdot
  R_3=0$. The remaining products involving~$e_6$,
  \begin{equation*}
    \begin{aligned}
      e_3e_6 &= -yz^4\e_{xz}\;,\\
      e_4e_6 &= -x^2yz^2\e_{xy} - z^5\e_{yz}\;,\text{ and}\\
      e_5e_6 &= -y^2z^3\e_{xz}
    \end{aligned}
  \end{equation*}
  are zero because the coefficients vanish in $R$. Thus, one has
  $A_1\cdot A_1 =0$.
\end{bfhpg*}

\begin{bfhpg*}[The product $A_1 \cdot A_2$]
  It follows from the computations in the previous section that the
  only non-zero products $\cls{e}_i\cls{f}_j$ for $1\le i \le 5$ and
  $j\in\set{1,2,4,6}$ are the ones listed in \eqref{I2G}.  The
  products $e_if_5$ for $1\le i \le 6$ are zero as one has $R_{\ge
    2}\cdot R_{4}=0$; similarly, the products
  \begin{equation*}
    e_1f_3\,,\; e_1f_7\,,\; e_2f_3\,,\; e_2f_7\,,\; e_6f_3\,,\; e_6f_4\,, \text{ and }\, 
    e_6f_7
  \end{equation*}
  are zero as one has $R_{3}\cdot R_{\ge 3}=0$. Among the remaining
  products, $e_4f_7$ and $e_6f_2$ are zero by graded-commutativity,
  and the following are zero because the coefficients vanish in $R$.
  \begin{alignat*}{2}
    e_3f_7 &= yz^4\e_{xyz}&\qquad
    e_5f_7 &= y^2z^3\e_{xyz}\\
    e_4f_3 &= z^5\e_{xyz}&
    e_6f_1 &= yz^4\e_{xyz}\\
    e_5f_3 &= y^5\e_{xyz}
  \end{alignat*}
  The one remainging product is $e_3f_3 = y^4z\e_{xyz} = g_1$.
\end{bfhpg*}

In terms of the basis
$\cls{e}_1,\ldots,\cls{e}_6,\cls{f}_1,\ldots,\cls{f}_7,\cls{g}_1,\cls{g}_2$
for $A_{\ge 1}$ the only non-zero products of basis vectors are
\begin{equation*}
  \cls{e}_i\cls{f}_i = \cls{g}_1 \quad \text{ for } 1\le i \le 3\;,
\end{equation*}
whence $R$ belongs to $\clG{3}$.

\section{Proof that $Q/\mathfrak{g}_3$ and $Q/\mathfrak{g}_4$ are type
  $2$ algebras in $\clG{4}$ and $\clG{5}$}
\label{sec:I45}

\noindent The arguments that show that $Q/\mathfrak{g}_3$ and
$Q/\mathfrak{g}_4$ are \textbf{G} algebras follow the argument in
\secref{I3} closely; we summarize them below.

\begin{bfhpg}[The quotient by $\mathfrak{g}_3$]
  \label{I4}
  Set $R=Q/\mathfrak{g}_3$; as one has
  $\mathfrak{g}_3=\mathfrak{g}_2+(x^2z^2)$ it follows from \eqref{I3}
  that the elements listed in the second column below form bases for
  the subspaces $R_n$.
  \begin{equation*}
    \text{\begin{tabular}{lll}
        $R_0$ & $1$\\
        $R_1$ & $x,\, y,\, z$\\
        $R_2$ & $x^2,\, xy,\, xz,\, y^2,\, yz,\, z^2$\\
        $R_3$ & $x^3,\, x^2y,\, x^2z,\, xz^2,\, y^3,\, y^2z,\, z^3\quad$\\
        $R_4$ & $x^4,\, x^3y,\, x^3z,\, xz^3$& $y^4=xz^3,\,y^3z=x^4,\,z^4=x^3y$ \\
        $R_5$ & $x^4y$& $xz^4=y^4z=x^4y$\\
      \end{tabular}}
  \end{equation*}
  It is straightforward to verify that the socle of $R$ is generated
  by the elements $x^4y$ and $x^3z$, so it has rank $2$. Set
  $A=\H{\Koszul}$; the homology classes of cycles
  \begin{equation*}
    g_1 = x^4y\e_{xyz} \qand g_2=x^3z\e_{xyz}
  \end{equation*}
  form a basis for $A_3$. One readily verifies that the homology
  classes of $e_1,\ldots,e_6$ and $f_1,f_2,f_3,f_5,f_7$ from
  \secref{I3}, see also \pgref{notation}, together with those of the
  cycles
  \begin{equation*}
    e_7 = xz^2\e_x\,, \quad f_4= y^3\e_{xy}-xz^2\e_{xz}\,, 
    \quad f_6=xz^2\e_{xy}\,,\qand\, f_8=x^3z\e_{xz} 
  \end{equation*}
  form bases for $A_1$ and $A_2$.

  The products $e_1e_7$, $e_2e_7$, and $e_6e_7$ vanish as one has
  $R_{3}\cdot R_{3} = 0$, while $e_3e_7$ and $e_5e_7$ are zero by
  graded-commutativity. Finally, one has $e_4e_7 = -xz^4\e_{xy} =
  -\partial(xz^3\e_{xyz})$, so also the product $\cls{e}_4\cls{e}_7$
  is zero. Together with the computations from the previous section
  this shows that $A_{1}\cdot A_{1}$ is zero.

  The products $e_7f_j$ vanish for $3\le j\le 8$ as one has
  $R_{3}\cdot R_{\ge 3} = 0$, and for the same reason any one of the
  elements $e_1$, $e_2$, and $e_6$ multiplied by either $f_4$ or $f_6$
  is zero. All products $e_if_8$ vanish as one has $R_{\ge 2}\cdot
  R_{4} = 0$. Among the remaining products involving $e_7$, $f_4$, or
  $f_6$ all but one are zero by graded-commutativity, the non-zero
  product is $e_4f_4 = xz^4\e_{xyz} = g_1$. It follows that in terms
  of the $\kk$-basis $\cls{e}_1,\ldots,\cls{e}_7$,
  $\cls{f}_1,\ldots,\cls{f}_8$, $\cls{g}_1,\cls{g}_2$ for $A_{\ge 1}$
  the only non-zero products of basis-elements are
  \begin{equation*}
    \cls{e}_i\cls{f}_i = \cls{g}_1 \quad \text{ for } 1\le i \le 4\;,
  \end{equation*}
  whence $R$ belongs to $\clG{4}$.
\end{bfhpg}

\begin{bfhpg}[The quotient by $\mathfrak{g}_4$]
  \label{I5}
  Set $R=Q/\mathfrak{g}_4$; as one has
  $\mathfrak{g}_4=\mathfrak{g}_3+(x^3z)$ it follows from \pgref{I4}
  that the elements listed in the second column below form bases for
  the subspaces $R_n$.
  \begin{equation*}
    \text{\begin{tabular}{lll}
        $R_0$ & $1$\\
        $R_1$ & $x,\, y,\, z$\\
        $R_2$ & $x^2,\, xy,\, xz,\, y^2,\, yz,\, z^2$\\
        $R_3$ & $x^3,\, x^2y,\, x^2z,\, xz^2,\, y^3,\, y^2z,\, z^3\quad$\\
        $R_4$ & $x^4,\, x^3y,\, xz^3$& $y^4=xz^3,\,y^3z=x^4,\,z^4=x^3y$ \\
        $R_5$ & $x^4y$& $xz^4=y^4z=x^4y$\\
      \end{tabular}}
  \end{equation*}
  It is straightforward to check that the socle of $R$ is generated by
  the elements $x^4y$ and $x^2z$, so it has rank $2$. Set
  $A=\H{\Koszul}$; the homology classes of the cycles
  \begin{equation*}
    g_1 = x^4y\e_{xyz} \qand g_2=x^2z\e_{xyz}
  \end{equation*}
  form a basis for $A_3$. One readily verifies that the homology
  classes of $e_1,\ldots,e_7$ and $f_1,f_2,f_3,f_4,f_6,f_7$ from
  \pgref{I4}, see also \pgref{notation}, together with those of
  \begin{equation*}
    e_8 = x^2z\e_x\,, \quad f_5= -x^3\e_{xz}+y^2z\e_{yz}\,, 
    \quad f_8=x^2z\e_{xy}\,,\qand\, f_9=x^2z\e_{xz} 
  \end{equation*}
  form bases for $A_1$ and $A_2$.

  All the products $e_ie_8$ are zero as $x^2z$ is in the socle of
  $R$. Together with the computations from \pgref{I4} this shows that
  $A_{1}\cdot A_{1}$ is zero.

  All products $e_8f_j$, $e_if_8$, and $e_if_9$ vanish as $x^2z$ is in
  the socle of $R$. Any one of the elements $e_1$, $e_2$, $e_6$, and
  $e_7$ multiplied by $f_5$ is zero for as one has $R_{3}\cdot R_{3} =
  0$. The remaining products involving $f_5$ are
  \begin{equation*}
    e_3f_5 = y^3z^2\e_{xyz} = 0\,,\quad
    e_4f_5 = x^3z^2\e_{xyz} = 0\,,\qand\,
    e_5f_5 = y^4z\e_{xyz} = g_1\;.
  \end{equation*}
  It follows that in terms of the $\kk$-basis
  $\cls{e}_1,\ldots,\cls{e}_8$, $\cls{f}_1,\ldots,\cls{f}_9$,
  $\cls{g}_1,\cls{g}_2$ for $A_{\ge 1}$ the only non-zero products of
  basis vectors are
  \begin{equation*}
    \cls{e}_i\cls{f}_i = \cls{g}_1 \quad \text{ for } 1\le i \le 5\;,
  \end{equation*}
  whence $R$ belongs to $\clG{5}$.
\end{bfhpg}

\section{Proof that $Q/\mathfrak{b}_1$ is a type $1$ algebra in
  $\clB$} 
\label{sec:J1}

\noindent The ideal $\mathfrak{b}_1=(x^3,x^2y,yz^2,z^3)$ in
$Q=\pows[\kk]{x,y,z}$ is generated by homogeneous elements, so
$R=Q/\mathfrak{b}_1$ is a graded $\kk$-algebra.  It is straightforward
to verify that the elements listed below form bases for the subspaces
$R_n$.
\begin{equation}
  \label{eq:J1}
  \text{\begin{tabular}{ll}
      $R_0$ & $1$\\
      $R_1$ & $x,\,y,\,z$\\
      $R_2$ & $x^2,\, xy,\, xz,\, y^2,\, yz,\, z^2$\\
      $R_3$ & $x^2z,\,xy^2,\,xyz,\,xz^2,\,y^3,\,y^2z$\\
      $R_4$ & $x^2z^2,\,xy^3,\,xy^2z,\,y^4,\,y^3z$\\
      $R_{n\ge 5}$ & $xy^{n-1},\,xy^{n-2}z,\,y^n,\,y^{n-1}z$
    \end{tabular}}
\end{equation}
Set $A=\H{\Koszul}$; we shall verify that $A$ has the multiplicative
structure described in \eqref{B} and that $R$ has socle rank $1$.

\begin{bfhpg*}[A basis for $A_{\ge 1}$]
  From \eqref{J1} it is straightforward to verify that the socle of
  $R$ is generated by $x^2z^2$, so it has rank $1$ and the homology
  class of the cycle
  \begin{equation}
    \label{eq:J1g}
    g_1 = x^2z^2\e_{xyz}
  \end{equation}
  is a basis for $A_3$.  The ideal $\mathfrak{b}_1$ is minimally
  generated by $4$ elements, whence one has $\rnk[\kk]{A_1}=4 =
  \rnk[\kk]{A_2}$; see \eqref{rank}. Proceeding as in \secref{I2} it
  is straightforward to verify that the elements $e_1,\ldots,e_4$ and
  $f_1,\ldots,f_4$ listed below are cycles in $\Koszul_1$ and
  $\Koszul_2$ whose homology classes form bases for the $\kk$-vector
  spaces $A_1$ and $A_2$.
  \begin{equation}
    \label{eq:J1ef}
    \begin{alignedat}{2}
      e_1 &= x^2\e_x&\qquad\qquad
      f_1 &= z^2\e_{yz}\\
      e_2 &= z^2\e_z&
      f_2 &= x^2\e_{xy}\\
      e_3 &= xy\e_x&
      f_3 &= x^2z^2\e_{xz}\\
      e_4 &= z^2\e_y& f_4 &= xyz\e_{xz}
    \end{alignedat}
  \end{equation}
\end{bfhpg*}

\begin{bfhpg*}[The product $A_1 \cdot A_1$]
  To determine the multiplication table $A_1\times A_1$ it is
  sufficient to compute the products $e_ie_j$ for $1\le i<j\le 4$; see
  \eqref{grcom}.  The product $e_1e_3$ is zero by graded-commutativity
  of $\Koszul$, and the following products are zero because the
  coefficients vanish in $R$.
  \begin{equation*}
    \begin{aligned}
      e_2e_3 &= -xyz^2\e_{xz}\\
      e_2e_4 &= -z^4\e_{yz}\\
      e_3e_4 &= xyz^2\e_{xy}
    \end{aligned}
  \end{equation*}
  The remaining products are
  \begin{equation*}
    e_1e_2 = x^2z^2\e_{xz} = f_3\qqand
    e_1e_4 = x^2z^2\e_{xy} = \partial(x^2z\e_{xyz})\;.
  \end{equation*}
  It follows that the product $\cls{e}_1\cls{e_4}$ in homology is
  zero, leaving us with a single non-zero product of basis vectors in
  $A_1$, namely $\cls{e}_1\cls{e}_2 = \cls{f}_3$.
\end{bfhpg*}

\begin{bfhpg*}[The product $A_1 \cdot A_2$]
  One has
  \begin{equation*}
    \begin{alignedat}{2}
      e_1f_1 &= x^2z^2\e_{xyz} = g_1&\qquad
      e_4f_3 &= -x^2z^4\e_{xyz} = 0\\
      e_2f_2 &= x^2z^2\e_{xyz} = g_1&
      e_4f_4 &= -xyz^3\e_{xyz} = 0\\
      &&e_3f_1 &= xyz^2\e_{xyz} = 0\;.
    \end{alignedat}
  \end{equation*}
  The remaining products are zero by graded-commutativity.
\end{bfhpg*}

In terms of the $\kk$-basis $\cls{e}_1,\ldots,\cls{e}_4,
\cls{f}_1,\ldots,\cls{f}_4,\cls{g}_1$ for $A_{\ge 1}$ the only
non-zero products of basis vectors are
\begin{equation}
  \label{eq:J1B}
  \begin{gathered}
    \cls{e}_1\cls{e}_2 =\cls{f}_3
  \end{gathered}
  \qqand
  \begin{gathered}
    \cls{e}_1\cls{f}_1 =\cls{g}_1\\
    \cls{e}_2\cls{f}_2 =\cls{g}_1\rlap{\;.}
  \end{gathered}
\end{equation}
It follows that $R$ belongs to the class $\clB$.

\section{Proof that $Q/\mathfrak{b}_2$ is a type $1$ algebra in
  $\clB$} 
\label{sec:J2}

\noindent Set $R=Q/\mathfrak{b}_2$; as one has
$\mathfrak{b}_2=\mathfrak{b}_1+(xyz)$ it follows from \eqref{J1} that
the elements listed below form bases for the subspaces $R_n$.
\begin{equation}
  \label{eq:J2}
  \text{\begin{tabular}{ll}
      $R_0$ & $1$\\
      $R_1$ & $x,\,y,\,z$\\
      $R_2$ & $x^2,\, xy,\, xz,\, y^2,\, yz,\, z^2$\\
      $R_3$ & $x^2z,\,xy^2,\,xz^2,\,y^3,\,y^2z$\\
      $R_4$ & $x^2z^2,\,xy^3,\,y^4,\,y^3z$\\
      $R_{n\ge 5}$ & $xy^{n-1},\,y^n,\,y^{n-1}z$
    \end{tabular}}
\end{equation}
Set $A=\H{\Koszul}$; we shall verify that $A$ has the multiplicative
structure described in \eqref{B} and that $R$ has socle rank $1$.

\begin{bfhpg*}[A basis for $A_{\ge 1}$]
  From \eqref{J2} it is straightforward to verify that the socle of
  $R$ is generated by $x^2z^2$, so it has rank $1$, and the homology
  class of the cycle $g_1$ from \eqref{J1g} is a basis for $A_3$;
  cf.~\pgref{notation}.  The ideal $\mathfrak{b}_2$ is minimally
  generated by $5$ elements, so one has $\rnk[\kk]{A_1}=5 =
  \rnk[\kk]{A_2}$; see \eqref{rank}. Proceeding as in \secref{I2} it
  is straightforward to verify that the homology classes of
  $e_1,\ldots,e_4$ from \eqref{J1ef} together with the class of the
  cycle
  \begin{equation*}
    e_5 = yz\e_x
  \end{equation*}
  form a basis for $A_1$. Similarly, one verifies that the homology
  classes of $f_1,f_2,f_3$ from \eqref{J1ef} together with those of
  the cycles
  \begin{equation*}
    f_4 = xy\e_{xz}\qqand f_5 = yz\e_{xz}
  \end{equation*}
  make up a basis for $A_2$.
\end{bfhpg*}

\begin{bfhpg*}[The product $A_1 \cdot A_1$]
  It follows from \eqref{J1B} that $\cls{e}_1\cls{e}_2=\cls{f}_3$ is
  the only non-zero product $\cls{e}_i\cls{e}_j$ for $1\le i < j \le
  4$. To complete the multiplication table $A_1\times A_1$ it is by
  \eqref{grcom} sufficient to compute the products $e_ie_5$ for $1\le
  i \le 4$. The products $e_1e_5$ and $e_3e_5$ are zero by
  graded-commutativity of $\Koszul$, and the remaining products
  involving $e_5$,
  \begin{equation*}
    e_2e_5 = -yz^3\e_{xz} \qqand e_4e_5 = -yz^3\e_{xy}\;,
  \end{equation*}
  are zero because the coefficients vanish in $R$.
\end{bfhpg*}

\begin{bfhpg*}[The product $A_1 \cdot A_2$]
  It follows from \eqref{J1B} that the only non-zero products
  $\cls{e}_i\cls{f}_j$ for $1 \le i \le 4$ and $1 \le j \le 3$ are
  $\cls{e}_1\cls{f}_1 = \cls{e}_2\cls{f}_2 = \cls{g}_1$. To complete
  the multiplication table $A_1 \times A_2$ one has to compute the
  products $e_5f_j$ for $1\le j \le 5$ and $e_if_4$ and $e_if_5$ for
  $1 \le i \le 5$. The next products are zero because the coefficients
  vanish in $R$,
  \begin{align*}
    e_4f_4 &= -xyz^2\e_{xyz}\\
    e_4f_5 &= -yz^3\e_{xyz}\\
    e_5f_1 &= yz^3\e_{xyz}\;,
  \end{align*}
  and the remaining are zero by graded-commutativity.
\end{bfhpg*}

In terms of the $\kk$-basis $\cls{e}_1,\ldots,\cls{e}_5,
\cls{f}_1,\ldots,\cls{f}_5,\cls{g}_1$ for $A_{\ge 1}$ the only
non-zero products of basis vectors are the ones listed in \eqref{J1B},
so $R$ belongs to $\clB$.

\section{Proof that $Q/\mathfrak{b}_3$ and $Q/\mathfrak{b}_4$ are type
  $3$ algebras in $\clB$}
\label{sec:J34}

\noindent The arguments that show that $Q/\mathfrak{b}_3$ and
$Q/\mathfrak{b}_4$ are $\clB$ algebras follow the argument in
\secref{J2} closely; we summarize them below.

\begin{bfhpg}[The quotient by $\mathfrak{b}_3$]
  \label{J3}
  Set $R=Q/\mathfrak{b}_3$; as one has
  $\mathfrak{b}_3=\mathfrak{b}_2+(xy^2-y^3)$ it follows from
  \eqref{J2} that the elements listed below form bases for the
  subspaces $R_n$.
  \begin{equation*}
    \text{\begin{tabular}{lll}
        $R_0$ & $1$\\
        $R_1$ & $x,\,y,\,z$\\
        $R_2$ & $x^2,\, xy,\, xz,\, y^2,\, yz,\, z^2$\\
        $R_3$ & $x^2z,\,xy^2,\,xz^2,\,y^2z$\qquad& $y^3=xy^2$\\
        $R_4$ & $x^2z^2$\\
      \end{tabular}}
  \end{equation*}
  It is straightforward to verify that the socle of $R$ is generated
  by the elements $x^2z^2$, $xy^2$, and $y^2z$, so it has rank
  $3$. Set $A=\H{\Koszul}$; the homology classes of the cycles
  \begin{equation*}
    g_1 = x^2z^2\e_{xyz}\,, \quad g_2 = xy^2\e_{xyz}\,, \qand\, g_3=y^2z\e_{xyz}
  \end{equation*}
  form a basis for $A_3$. One readily verifies that the homology
  classes of $e_1,\ldots,e_5$ and $f_1,\ldots,f_5$ from the previous
  section, see also \pgref{notation}, together with those of the
  cycles
  \begin{equation*}
    e_6 = y^2\e_x-y^2\e_y\,, \quad f_6=xy^2\e_{xy}\,, 
    \quad f_7=y^2z\e_{xy}\,,\qand\, f_8=y^2z\e_{yz} 
  \end{equation*}
  form bases for $A_1$ and $A_2$.

  The products $e_ie_6$, $e_if_6$, $e_if_7$, and $e_if_8$ for $1\le
  i\le 5$ as well as the products $e_6f_j$ for $1 \le j \le 8$ vanish
  as one has $y^2R_{\ge2}=0$. It follows that in terms of the
  $\kk$-basis $\cls{e}_1,\ldots,\cls{e}_6,
  \cls{f}_1,\ldots,\cls{f}_8,\cls{g}_1,\cls{g}_2,\cls{g}_3$ for
  $A_{\ge 1}$ the only non-zero products of basis vectors are the ones
  listed in \eqref{J1B}, so $R$ belongs to $\clB$.
\end{bfhpg}

\begin{bfhpg}[The quotient by $\mathfrak{b}_4$]
  \label{J4}
  Set $R=Q/\mathfrak{b}_4$; as one has
  $\mathfrak{b}_4=\mathfrak{b}_3+(y^2z)$ it follows from \pgref{J3}
  that the elements listed below form bases for the subspaces $R_n$.
  \begin{equation*}
    \text{\begin{tabular}{lll}
        $R_0$ & $1$\\
        $R_1$ & $x,\,y,\,z$\\
        $R_2$ & $x^2,\, xy,\, xz,\, y^2,\, yz,\, z^2$\\
        $R_3$ & $x^2z,\,xy^2,\,xz^2$\qquad& $y^3=xy^2$\\
        $R_4$ & $x^2z^2$\\
      \end{tabular}}
  \end{equation*}
  It is straightforward to check that the socle of $R$ is generated by
  the elements $x^2z^2$, $xy^2$, and $yz$, so it has rank $3$. Set
  $A=\H{\Koszul}$; the homology classes of the cycles
  \begin{equation*}
    g_1 = x^2z^2\e_{xyz}\,, \quad g_2 = xy^2\e_{xyz}\,, \qand\, g_3=yz\e_{xyz}
  \end{equation*}
  form a basis for $A_3$. One readily verifies that the homology
  classes of $e_1,\ldots,e_6$ and $f_1,\ldots,f_6$ from \pgref{J3}
  together with those of
  \begin{equation*}
    e_7 = yz\e_y\, \quad f_7=yz\e_{xy}\,, 
    \quad f_8=yz\e_{yz}\,,\qand \, f_9=y^2\e_{xz}-y^2\e_{yz} 
  \end{equation*}
  form bases for $A_1$ and $A_2$.

  The products $e_ie_7$ for $1\le i\le 6$ and $e_7f_j$ for $1 \le j
  \le 9$ vanish as $yz$ is in the socle of $R$, and for $1\le i\le 6$
  the products $e_if_7$ and $e_if_8$ vanish for the same
  reason. Finally, all products $e_if_9$ vanish as one has $y^2R_{\ge
    2}=0$. Thus, in terms of the $\kk$-basis
  $\cls{e}_1,\ldots,\cls{e}_7,
  \cls{f}_1,\ldots,\cls{f}_9,\cls{g}_1,\cls{g}_2,\cls{g}_3$ for
  $A_{\ge 1}$ the only non-zero products of basis vectors are the ones
  listed in \eqref{J1B}, so $R$ belongs to $\clB$.
\end{bfhpg}

\section*{Acknowledgments} 

\noindent
We thank Jerzy Weyman for insightful discussions and suggestions and
we thank R\u{a}zvan Veliche for his help during the writing phase of
the project. The computer algebra software MACAULAY 2 \cite{M2} was
used to conduct the experiments that led to the discovery of the
algebras in Theorems I and II.

\bibliographystyle{amsplain}


\def\cprime{$'$}
  \providecommand{\arxiv}[2][AC]{\mbox{\href{http://arxiv.org/abs/#2}{\sf
  arXiv:#2 [math.#1]}}}
  \providecommand{\oldarxiv}[2][AC]{\mbox{\href{http://arxiv.org/abs/math/#2}{\sf
  arXiv:math/#2
  [math.#1]}}}\providecommand{\MR}[1]{\mbox{\href{http://www.ams.org/mathscinet-getitem?mr=#1}{#1}}}
  \renewcommand{\MR}[1]{\mbox{\href{http://www.ams.org/mathscinet-getitem?mr=#1}{#1}}}
\providecommand{\bysame}{\leavevmode\hbox to3em{\hrulefill}\thinspace}
\providecommand{\MR}{\relax\ifhmode\unskip\space\fi MR }
\providecommand{\MRhref}[2]{%
  \href{http://www.ams.org/mathscinet-getitem?mr=#1}{#2}
}
\providecommand{\href}[2]{#2}

\end{document}